\documentclass[12pt]{article}
\usepackage{amssymb,amsmath}
\usepackage[dvips]{color}
\usepackage[dvips]{graphics}
\usepackage{graphicx}

\newtheorem{theorem}{Theorem}[section]
\newtheorem{lemma}{Lemma}[section]

\newcommand{\numlist}[3]{#1=#2, \ldots , #3}

\newcommand{\normF}[1]{\|#1\|_F}

\newcommand{\ysqrt}[1]{\left( #1 \right)^{1/2}}
\newcommand{\ysqrtn}[1]{\left(#1 \right)^{-1/2}}

\newcommand{\normtwo}[1]{\|#1\|_2}
\newcommand{\normal}[1]{#1^T\!~#1}

\newcommand{\XXT}[1]{#1~#1^T}
\newcommand{\macheps}{\varepsilon_M}

\newcommand{\hot}{+ {\mathcal O(\macheps^2 )}}

\newcommand{\Real}{\mathbb{R}}
\newcommand{\range}{\mathrm{Range}}
\newcommand{\mbf}{\mathbf}     
\newcommand{\tbf}{\textbf}

\begin{document}
\title{Reorthogonalized Block Classical Gram--Schmidt}
\author{Jesse L. Barlow \thanks{The research of Jesse L. Barlow was sponsored by the National Science Foundation under contract no. CCF-1115704.}
\thanks{Department of Computer Science and Engineering, University Park, PA 16802-6822, USA, e-mail: barlow@cse.psu.edu.}    \and
Alicja Smoktunowicz \thanks{Faculty of Mathematics and Information Science,  Warsaw University of Technology,  Pl. Politechniki 1, 00-661,  Warsaw, Poland , e-mail: smok@mini.pw.edu.pl.}}
\maketitle

\begin{abstract}
A new reorthogonalized block classical Gram--Schmidt algorithm  is proposed that factorizes a full column rank matrix $A$ into $A=QR$ where $Q$ is left orthogonal (has orthonormal columns) and $R$ is upper triangular and nonsingular.

With appropriate assumptions on the diagonal blocks of $R$, the algorithm, when implemented in
floating point arithmetic with machine unit $\macheps$, produces $Q$ and $R$ such that
$\| I- Q^{T} Q \|_2 =O(\macheps)$ and $\| A-QR \|_2 =O(\macheps \| A \|_2)$. The resulting bounds also improve a previous bound by Giraud et al. [Num. Math., 101(1):87-100,\ 2005] on the CGS2 algorithm originally developed by Abdelmalek [BIT, 11(4):354--367,\ 1971].

\medskip

\end{abstract}

\section{Introduction}
\label{intro} \noindent
For a matrix $A \in  \mathbb{R}^{m\times n}$, $m \geq n$, we consider the computation of the Q--R decomposition
\begin{equation}
A = QR
\label{eq:QRdecomp}
\end{equation}
where $Q \in  \mathbb{R}^{m\times n}$ is left orthogonal (i.e., $Q^{T} Q=I_n$) and $R \in  \mathbb{R}^{n\times n}$ is upper triangular.
The matrix $A$ is assumed to have full column rank.

The approach considered is block classical Gram--Schmidt with reorthogonalization (BCGS2)
which operates on groups of columns of $A$ instead of columns in order to create a BLAS--3 \cite{DoDu90} compatable algorithm.
Thus we assume that $A$ is partitioned into blocks $A=(A_1, \ldots, A_s)$, where each $A_i$ has $p_i$ columns, i.e. $A_i(m \times p_i)$ for
$i=1, \ldots, s$, with $n=p_1+p_2+ \ldots + p_s$.

The block Gram--Schmidt algorithm that we present is a generalization of the classical Gram--Schmidt method with reorthogonalization (CGS2) which was first proposed and analyzed by Abdelmalek \cite{Abd71}, but our analysis and development follow the flavor of that given by Giraud et al \cite{Gir05}. A similar block algorithm based upon CGS, justified only by numerical tests, is proposed by Stewart \cite{Stew08}.
Other block Gram--Schmidt algorithms are presented by Jalby and Phillippe \cite{JaPh91} and Vanderstaeten \cite{Van00}.
 A goal of this paper is to establish the stability analysis results for BCGS2.
As summarized in \S \ref{sec:CGS2}, our analysis has implications for the algorithm in \cite{Abd71,Gir05} and shows that the CGS2 algorithms produces a near orthogonal matrix under weaker assumptions than given in \cite{Gir05}. An excellent summary of the role of Gram--Schmidt algorithms is
given in \cite[\S 2.4, \S 3.2]{Bjor96} and in the 1994 survey  paper \cite{Bjor94}.

We prove that  BCGS2 is numerically stable under natural conditions outlined in \S \ref{sec:CGS2erranalysis} producing computed $Q$ and $R$ by BCGS2 in floating point arithmetic with machine unit $\macheps$ that satisfy
\begin{eqnarray}
\| I_n -Q^{T} Q \|_2 &\leq& \macheps f_1(m,n,p)  \hot, \label{eq:IQTQbnd}\\
\| A-QR \|_2 &\leq&  \macheps f_2(m,n,p)  \| A \|_2 \hot \label{eq:AQRbnd}
\end{eqnarray}
for modestly growing functions $f_1(\cdot)$ and $f_2(\cdot)$,  where $p = \max_{1\leq i\leq s} p_i$  is the maximum block size of $A$. We assume that $f_j(m,n,p)\macheps <1, j=1,2$, otherwise these bounds are meaningless.

At the core of our block CGS method is a routine \tbf{local\_qr}, where for a matrix $B \in \mathbb{R}^{m\times p}$, $p \leq n \leq m$, produces
\[
[Q,R] =\mbf{local\_qr}(B)
\]
with
\begin{eqnarray}
\| I_p-Q^{T} Q \|_2  &\leq& \macheps L_1(m,p) <1 ,  \label{eq:IQBbnd} \\
\| B-Q R \|_2  &\leq&  \macheps L_1(m,p)  \| B \|_2  \label{eq:BQRbnd}
\end{eqnarray}
for some modest function $L_1(\cdot)$.
The routine \tbf{local\_qr} may be produced using Householder or Givens Q--R factorization.
For appropriate BLAS--3 speed \cite{DoDu90}, that is, to take advantage of caching, the implementation of \tbf{local\_qr} may be done using the ``tall, skinny'' Q--R (TSQR) discussed in the recent Ph.D. thesis by  Hoemmen \cite[\S 2.3]{HoePhD10}.
An interpretaton of \cite[\S 19.3]{Hig02} on the error analysis of Householder Q--R yield a
function $L_1(m,p)=d_1 mp^{3/2}$ where $d_1$ is a constant.
If we make the assumptions in Theorem \ref{thm:IQTQbnd} for $p=1$, then the CGS2 algorithm
described by Giraud et al. \cite{Gir05} satisfies (\ref{eq:IQBbnd})--(\ref{eq:BQRbnd}).

In \S \ref{sec:algorithm}, we present our algorithm, and in \S \ref{sec:CGS2erranalysis}, we
prove the properties (\ref{eq:IQTQbnd})-- (\ref{eq:AQRbnd}), followed by our conclusions in \S \ref{sec:conclusions}.

\section{Reorthogonalized Block Gram--Schmidt}
\label{sec:algorithm}

First, we summarize the algorithm in \cite{Abd71,Gir05}. Let $A \in  \mathbb{R}^{m\times n}$ be given column-wise as
\[
A=  (\mbf{a}_1, \ldots , \mbf{a}_n).
\]

To generate the decomposition (\ref{eq:QRdecomp}) with $R=(r_{ij})$ and $Q = (\mbf{q}_1, \ldots , \mbf{q}_n)$, a step in the CGS2 algorithm takes a near left orthogonal matrix $U \in \mathbb{R}^{m\times p}$,  a
vector $\mbf{b} \in \Real^m$, and produces $r_b \in \Real$, $\mbf{s}_b \in \Real^{p}$ and $\mbf{q}_b\in \Real^m$ such that
\begin{eqnarray}
\mbf{q}_b r_b &=& (I_m-\XXT{U})^2\mbf{b}, \quad \|{\mbf{q}_b}\|_2 =1, \label{eq:CGS2eq}\\
\mbf{b} &=& U\mbf{s}_b + \mbf{q}_b r_b. \label{eq:lineareq}
\end{eqnarray}
Ideally, $\mbf{q}_b$ is orthogonal to the columns of $U$, so  (\ref{eq:CGS2eq}) should be replaced by
\begin{equation}
U^T \mbf{q}_b =0, \quad \|{\mbf{q}_b}\|_2 =1. \label{eq:orthrel}
\end{equation}
However, the condition (\ref{eq:orthrel}) becomes hard to enforce as $r_b$ approaches zero. Unlike the version of this procedure in \cite{Gir05}, we scale the approximation to $\mbf{q}_b$ at each step. That change has two benefits: (1)
the function \tbf{cgs2\_step} below is more resistant to underflow; (2) it leads to a natural generalization to a block algorithm.

\medskip

\noindent
{\bf Function 2.1 [One step of CGS2]} $\:$ \label{fun:CGS2_step}

\begin{tabbing}
\tbf{function} $[\mbf{q}_b,r_b,\mbf{s}_b]$ = {\em cgs2\_step}($U,\mbf{b}$) \\
$\mbf{s}_1 = U^T \mbf{b}$; $\mbf{y}_1 = \mbf{b}-U\mbf{s}_1$; \\
$r_1 = \|{\mbf{y}_1}\|_2$; $\mbf{q}_1 = \mbf{y}_1/r_1$;\\
$\mbf{s}_2 = U^T \mbf{q}_1$; $\mbf{y}_2 = \mbf{q}_1 - U\mbf{s}_2$;\\
$r_2 = \|{\mbf{y}_2}\|_2$; $\mbf{q}_b = \mbf{y}_2/r_2$;\\
$\mbf{s}_b = \mbf{s}_1 + \mbf{s}_2 r_1$; $r_b = r_2 r_1$; \\
\tbf{end} {\em cgs2\_step} \\
\end{tabbing}

Notice that one step of CGS2 consists of exactly two steps of CGS. We have
\[
\mbf{q}_1 r_1 = (I_m-\XXT{U}) \mbf{b}, \quad  \mbf{q}_b r_2 = (I_m-\XXT{U}) \mbf{q_1},
\]
so, clearly,  (\ref{eq:CGS2eq}) holds.

The CGS2 algorithm from \cite{Abd71,Gir05} for computing the  Q--R decomposition is stated next.

\medskip

\noindent
{\bf Function 2.2 [Classical Gram--Schmidt with Reorthogonalization (CGS2)]}  $\:$ \label{fun:CGS2}

\begin{tabbing}
\tbf{function}  $[Q,R]$={\em cgs2}($A$) \\
$[m,n]$=\tbf{size}($A$); \\
$R=\|{A(:,1)}\|_2$; $Q=A(:,1)/R$; \\
\tbf{for} \= $k=2:n$ \\
\>$[\mbf{q}_{new},r_{new},\mbf{s}_{new}]$=\tbf{cgs2\_step}($Q,A(:,k)$);\\
\>$R = \left(\begin{array}{cc} R & \mbf{s}_{new} \\ 0 & r_{new} \end{array}\right)$; \,\,\, $Q = \left(\begin{array}{cc} Q & \mbf{q}_{new} \end{array}\right)$; \\
\tbf{end}; \\
\tbf{end}; {\em cgs2} \\
\end{tabbing}

To obtain the new function  \tbf{block\_CGS2\_step}, the block analog of \tbf{cgs2\_step}, the function \tbf{local\_qr} substitutes for scaling the vectors. First we introduce the function \tbf{block\_CGS\_step}.
Upon inputting $B \in  \mathbb{R}^{m\times p}$, $U \in \mathbb{R}^{m\times t}$,  $r +p \leq n \leq m$,  we produce
$\bar{Q} \in  \mathbb{R}^{m\times p}$, $\bar{R} \in \mathbb{R}^{p\times p}$, and $\bar{S} \in \mathbb{R}^{t\times p}$ such that
\begin{eqnarray}
\bar{Q} \bar{R} &=& (I_m-\XXT{U}) B, \quad  \normal{\bar{Q}} =I_p, \label{eq:blockCGSeq}\\
B &=& U\bar{S} + \bar{Q} \bar{R}. \label{eq:blocklineareq1}
\end{eqnarray}

\medskip

\noindent
{\bf Function 2.3 [One step of block CGS]} $\:$ \label{fun:BCGSstep}
\begin{tabbing}
\tbf{function} $[\bar{Q},\bar{R},\bar{S}]$ = {\em block\_CGS\_step}($U,B$) \\
$\bar{S} = U^T B$; \\
$\bar{Y} = B-U \bar{S}$;\\
$[\bar{Q},\bar{R}]$=\tbf{local\_qr}($\bar{Y}$); \\
\tbf{end} {\em block\_CGS\_step} \\
\end{tabbing}

Since we assume that $U$ and $\bar{Q}$ are near orthogonal in the sense that
\begin{eqnarray}
\normtwo{I_t-\normal{U}} &\leq& \macheps f_1(m,t,p) \hot<1, \label{eq:Uorth} \\
\normtwo{I_p-\normal{\bar{Q}}} &\leq& \macheps L_1(m,p)\hot, \label{eq:Qorth}
\end{eqnarray}
a simple eigenvalue/singular value analysis yields the bounds
\begin{eqnarray*}
\normtwo{U}&\leq& 1+0.5 \macheps f_1(m,t,p) \hot =1+\mathcal{O}(\macheps), \\
\normtwo{\bar{Q}} &\leq& 1+0.5\macheps L_1(m,p) \hot =1+\mathcal{O}(\macheps)
\end{eqnarray*}
of which we will make generous use  throughout our analysis.

The behavior of this routine in floating point arithmetic is given by the next two lemmas.
The proof of the first one is elementary, obvious, and will be skipped.

\begin{lemma}
\label{lem:BCGS_step} In floating point arithmetic with machine unit $\macheps$, Function 2.3
produces $\bar{Q}$,$\bar{R}$,$\bar{S}$ and $\bar{Y}$ such that for $L_1(\cdot)$ defined in \eqref{eq:IQBbnd}--\eqref{eq:BQRbnd}
and for modestly growing functions $L_2(\cdot)$ and $L_3(\cdot)$ we have
\begin{eqnarray}
\normtwo{I_p-\normal{\bar{Q}}} &\leq& \macheps L_1(m,p) \\
\bar{Q} \bar{R} &=& \bar{Y} + \Delta \bar{Y}, \quad \normtwo{\Delta \bar{Y}} \leq \macheps L_1(m,p) \normtwo{B} \\
\bar{S} + \delta \bar{S} &=& U^T B, \quad \normtwo{\delta \bar{S}} \leq \macheps L_2(m,t,p) \normtwo{B} \hot \\
\bar{Y} + \delta \bar{Y} &=& B-U\bar{S}, \quad \normtwo{\delta \bar{Y}} \leq \macheps L_3(t,p) \normtwo{B} \hot
\end{eqnarray}
\end{lemma}

If we use a standard matrix multiply and add routine, reasonable values for $L_2(m,t,p)$ and $L_3(t,p)$ are
\[
L_2(m,t,p)=m t^{1/2} p^{1/2}, \quad L_3(t,p)= p^{1/2}(1+t^{3/2}).
\]

A second lemma yields a backward error bound.

\begin{lemma}
\label{lem:BCGS_step_backward} Assume the hypothesis and notation of Lemma \ref{lem:BCGS_step}, then the computed
$\bar{Q}$ and $\bar{R}$ from Function 2.3
satisfy
\begin{equation}
\bar{Q} \bar{R} = (I_m-\XXT{U}) B +  F
\label{eq:Fdef1}
\end{equation}
where
\begin{equation}
  F  = \Delta \bar{Y} -\delta \bar{Y}+U(\delta \bar{S}).
\label{eq:Fdef2}
\end{equation}
Thus
\begin{equation}
\normtwo{F} \leq \macheps L_F(m,t,p) \normtwo{B} \hot,
\label{eq:Fbnd1}
\end{equation}
where
\begin{equation}
L_F(m,t,p)= L_1(m,p)+L_2(m,t,p)+L_3(t,p).
\label{eq:Fbnd2}
\end{equation}
\end{lemma}

\noindent

{\bf Proof.}
We simply unwind the relationships from Lemma \ref{lem:BCGS_step} to obtain
\begin{eqnarray*}
\bar{Q} \bar{R} &=& \bar{Y}+ \Delta \bar{Y} \\
&=& B-U\bar{S}-\delta \bar{Y}+ \Delta \bar{Y} \\
&=& (I_m-\XXT{U})B + U(\delta \bar{S})-\delta \bar{Y}+ \Delta \bar{Y} \\
&=& (I_m-\XXT{U})B +F
\end{eqnarray*}
which is (\ref{eq:Fdef1})--(\ref{eq:Fdef2}). The use of standard norm inequalities yields (\ref{eq:Fbnd1})--(\ref{eq:Fbnd2}).

\medskip

A norm bound that results from Lemma \ref{lem:BCGS_step_backward} is necessary for
our error analysis.

\begin{lemma}
\label{lem:barRbnd}
Assume the hypothesis and notation of Lemma \ref{lem:BCGS_step}, then the computed $\bar{R}$ from Function 2.3
satisfies
\begin{equation}
\normtwo{\bar{R}} \leq \normtwo{B}(1+ \mathcal{O}(\macheps)).
\label{eq:barRbnd}
\end{equation}
\end{lemma}

\noindent

{\bf Proof.}
Taking (\ref{eq:Fdef1}) multiplying on the left by $Q^T$ and reorganizing terms, we have
\[
 \bar{R} = \bar{Q}^T (I_m-\XXT{U}) B + (I_p-\normal{\bar{Q}})\bar{R}+ \bar{Q}^T F
\]
thus
\begin{equation}
\normtwo{\bar{R}} \leq \normtwo{I_m-\XXT{U}} \normtwo{B} + \normtwo{I_p-\normal{\bar{Q}}}\normtwo{\bar{R}} + \normtwo{F} \hot
\label{eq:Rbnd}
\end{equation}
To bound $\normtwo{\bar{R}}$ in (\ref{eq:Rbnd}),
we first need to bound $\normtwo{I_m-\XXT{U}}$. Let $U$ have the Q--R decomposition
\[
U=Z\left(\begin{array}{c} R_U \\ 0_{(m-t)\times t }\end{array}\right)
\]
where $Z$ is orthogonal and $R_U$ is upper triangular. Then
\[
\normtwo{I_t-\normal{U}} = \normtwo{I_t-\normal{R_U}} \leq \macheps f_1(m,t,p) \hot
\]
and
\[
I_m-\XXT{U} = \left(\begin{array}{cc} I_t -\XXT{R_U} & 0 \\ 0 & I_{m-t} \end{array}\right).
\]
Since $\normal{R_U}$ and $\XXT{R_U}$ have the same eigenvalues,
\[
\normtwo{I_t-\XXT{S}} \leq \macheps f_1(m,t,p)\hot.
\]
Using the assumption (\ref{eq:Uorth}) and the results of Lemma \ref{lem:BCGS_step_backward},
we have
\begin{eqnarray}
\normtwo{I_m-\XXT{U} } &=& \max\{ \normtwo{I_t-\XXT{R_U}} ,1\}\nonumber \\
 &=& \max\{ \normtwo{I_t-\normal{R_U}},1\}\nonumber \\
 &\leq& \max\{ \macheps f_1(m,t,p),1\} = 1 \label{eq:projbnd}
\end{eqnarray}
by our assumption about $f_1(m,t,p)$ in (\ref{eq:IQTQbnd}). Thus, excepting $\mathcal{O}(\macheps^2)$ terms,
\begin{equation}
\normtwo{\bar{R}} \leq  \normtwo{I_m-\XXT{U}}\normtwo{B} + \macheps L_F(m,t,p)\normtwo{B}+\macheps L_1(m,p) \normtwo{\bar{R}}.
\label{eq:barRnormineq}
\end{equation}
Using (\ref{eq:projbnd}), (\ref{eq:barRnormineq}), and solving for $\normtwo{\bar{R}}$ yields
\begin{eqnarray*}
\normtwo{\bar{R}} &\leq& [(1+\macheps L_F(m,t,p))/(1-\macheps L_1(m,p))] \normtwo{B} \nonumber \\
&=& (1+\macheps[L_F(m,t,p)+L_1(m,p)]\hot)\normtwo{B} \\
&=& \normtwo{B} (1+ \mathcal{O}(\macheps))
\end{eqnarray*}
which is (\ref{eq:barRbnd}).

\medskip

Now we introduce the function  \tbf{block\_CGS2\_step} which consists of two steps of \tbf{block\_CGS\_step}.
The first step  of \tbf{block\_CGS2\_step} for given $B$ produces $Q_1, R_1, S_1$ such that
\begin{equation}
B= U S_1 + Q_1 R_1.
\label{eq:Beq}
\end{equation}

The second step takes  $Q_1$ and yields $Q_B, R_2, S_2$ satisfying
\begin{equation}
Q_1= U S_2 + Q_B R_2.
\label{eq:Q1eq}
\end{equation}

From (\ref{eq:Beq})--(\ref{eq:Q1eq}), it follows that
\begin{equation}
B= U S_1 + (Q_B R_2 + U S_2)R_1 = U (S_1 + S_2 R_1) + Q_B (R_2 R_1),
\end{equation}
so
\begin{equation}\label{eq:blocklineareq}
B = US_B + Q_B R_B
\end{equation}
where
\[
S_B = S_1+S_2 R_1, \quad R_B = R_2 R_1.
\]

\medskip

\noindent
{\bf Function 2.4 [One step of block CGS2]} $\:$ \label{fun:BCGS2_step}
\begin{tabbing}
\tbf{function} $[Q_B,R_B,S_B]$ = {\em block\_CGS2\_step}($U,B$) \\
 $[Q_1,R_1,S_1]$ = {\em block\_CGS\_step}($U,B$);\\
 $[Q_B,R_2,S_2]$ = {\em block\_CGS\_step}($U,Q_1$);\\
$S_B = S_1 + S_2 R_1$;\\
$R_B = R_2 R_1$; \\
\tbf{end} {\em block\_CGS2\_step} \\
\end{tabbing}

In exact arithmetic,
\begin{equation}
Q_B R_B = (I_m-\XXT{U})^2 B,
 \quad  {Q_B}^{T} {Q_B} =I_p. \label{eq:blockCGS2eq}
\end{equation}

Thus, we expect that
\begin{eqnarray*}
\|{U^T Q_B}\|_2 &\leq& \|{I_t-U^{T} U}\|_2^2  \|{ U^T B R_B^{-1}}\|_2\\
&\leq& \| I_t-U^{T} U\|_2^2  \normtwo{U}         \| B\|_2  \| R_B^{-1}\|_2\\
&\leq& \| I_t-U^{T} U\|_2^2   \| B\|_2  \| R_B^{-1}\|_2.
\end{eqnarray*}

Thus, as in (\ref{eq:CGS2eq})--(\ref{eq:orthrel}), ideally, (\ref{eq:blockCGS2eq}) should be replaced by
\begin{equation}
U^T Q_B =0, \quad {Q_B}^{T} {Q_B} =I_p, \label{eq:blockorthrel}
\end{equation}
but (\ref{eq:blockorthrel}) is hard to enforce if $\|{U^T B R_{B}^{-1}}\|_2$ is too large. When accounting for
rounding error, we can only guarantee that $\normtwo{U^T Q_B}$ is small if we bound $\normtwo{B} \normtwo{R_B^{-1}}$ as
shown in (\ref{eq:assumption1}).

Also, we have two bounds that result from interpreting Lemma \ref{lem:barRbnd}. They
are
\begin{eqnarray}
\normtwo{R_1} &\leq& \normtwo{B}(1+\mathcal{O}(\macheps)) \label{eq:R1bnd} \\
\normtwo{R_2} &\leq& \normtwo{Q_1}(1+\mathcal{O}(\macheps)) \leq  1+\mathcal{O}(\macheps) \label{eq:R2bnd}
\end{eqnarray}
which are freely used in our analysis.

In exact arithmetic,
\[
\left(\begin{array}{cc} U & B\end{array}\right) =\left(\begin{array}{cc} U &  Q_B\end{array}\right)  \left(\begin{array}{cc} I_t  & S_{B} \\ 0 & R_{B} \end{array}\right).
\]
An approach to developing a procedure similar to Function 2.4
is given by Strathopoulos and Wu \cite{StWu02}. In our notation,
they find $Q_B$ such that
\[
\range[(\begin{array}{cc} U & B \end{array})] = \range[(\begin{array}{cc} U & Q_B\end{array})]
\]
that also satisfies (\ref{eq:blockorthrel}). The focus of their paper is a procedure for \tbf{local\_qr} that is designed to be efficient in terms of storage accesses
if $B$ is ``tall and thin'' (i.e., $m \gg p$), but satisfies neither of the criteria (\ref{eq:IQBbnd})--(\ref{eq:BQRbnd}) well. The authors compensate by crafting a routine like that above, but
 using outer iterations to get $\|{U^T Q_B}\|_2$ as small  as possible and inner iterations to get $\|I_p -{Q_B}^{T} {Q_B}\|_2$ as small as possible. For our routine, the number of inner iterations is $1$ and outer iterations is $2$. The concern about ``tall and thin''
matrices $B$ is alleviated by use of the ``tall, skinny'' Q--R (TSQR) as in \cite[\S 2.3]{HoePhD10}.

As well as repeating the two operations from \tbf{block\_CGS\_step},
 there are two other operations for which we need error bounds.
In floating point arithmetic, the computed values of $S_B$ and $R_B$ from Function 2.5
satisfy
\begin{eqnarray}
S_B + \delta S_B &=& S_1 + S_2 R_1, \quad \normtwo{\delta S_B} \leq \macheps L_4(p)\normtwo{B}, \label{eq:deltaSB} \\
R_B + \delta R_B &=& R_2 R_1, \quad \normtwo{\delta R_B} \leq \macheps L_5(p) \normtwo{B}, \label{eq:deltaRB}
\end{eqnarray}
where $L_4(\cdot)$ and $L_5(\cdot)$ are modestly growing functions. For conventional matrix multiply and add,
$L_4(p)= p^{1/2}(1+p^{3/2})$ and $L_5(p)=p^2$. Using Weyl's inequality for singular values \cite[Corollary 8.6.2]{GoVa96}, we have
\begin{equation}
|\sigma_{\ell} (R_B) - \sigma_{\ell}(R_2 R_1)| \leq \normtwo{\delta R_B} \leq \macheps L_5(p) \normtwo{B},\quad \numlist{\ell}{1}{p}.
\label{eq:singvalbnd}
\end{equation}

In the block analog of Function 2.2,
we partition $A \in  \mathbb{R}^{m\times n}$  into
\begin{equation}
A = (A_1, A_2 , \ldots, A_s)
\label{eq: Ablocks}
\end{equation}
where  $A_k \in  \mathbb{R}^{m\times p_k}$ for $k=1, \ldots, s$.  In practice, if $p=\lceil n/s \rceil$, then usually, $p_i \in \{ p-1,p\}$. In this input to Function 2.5, we define the
parameter $\mbf{blocks}$ as
\begin{equation}
\mbf{blocks}=(p_1, \ldots, p_s)^T.
\label{eq:blocksdef}
\end{equation}

Likewise, we partition $Q$ into
\begin{equation}
Q = (Q_1, Q_2 , \ldots, Q_s).
\label{eq:Qblocks}
\end{equation}

We let
\begin{equation}
\hat{Q}_k =  (Q_1, Q_2 , \ldots, Q_k), \quad  \hat{A}_k =  (A_1, A_2 , \ldots, A_k)
\label{eq:hatQk}
\end{equation}
and let
\begin{equation}
R_k =\left(\begin{array}{cccccc} R_{11} & R_{12} & \cdots & \cdots & \cdots &  R_{1k} \\ & R_{22} & \cdots & \cdots & \cdots & R_{2k} \\  && \cdots & \cdots & \cdots &\cdots \\
&&&&& R_{kk} \end{array}\right).
\label{eq:Rk}
\end{equation}

The initial step of factoring the first block is $[Q_1,R_{1}]$=\tbf{local\_qr}($A_1$) with $R_{11}=R_1$ and $\hat Q_1=Q_1$.
Then for $k=1, \ldots, s-1$ we compute $S_{k+1}$, $Q_{k+1}$ and $R_{k+1,k+1}$ by Function 2.4 for $B=A_{k+1}$ and $U=\hat{Q}_k$.
Thus if
\begin{equation}
R_{k+1}= \left(\begin{array}{cc} R_k & S_{k+1} \\ 0 & R_{k+1,k+1} \end{array}\right), \quad \hat{Q}_{k+1} = \left(\begin{array}{cc} \hat {Q}_{k}, & Q_{k+1} \end{array}\right), \quad
\hat {A}_{k+1} = \left(\begin{array}{cc} \hat {A}_{k} & A_{k+1} \end{array}\right),
\end{equation}
then
\[
 \hat A_{k+1}=\hat{Q}_{k+1} R_{k+1}, \numlist{k}{1}{s-1}
\]
and, finally, $A=QR$, with $Q=\hat{Q}_s$ and $R=R_s$.

We summarize the algorithm BCGS2 as follows.

\medskip

\noindent
{\bf Function 2.5 [Block Classical Gram--Schmidt with Reorthogonalization (BCGS2)]} $\:$ \label{fun:blockCGS2}

\begin{tabbing}
\tbf{function}  $[Q,R]$={\em block\_CGS2} ($A,\mbf{blocks}$) \\
$[m,n]$=\tbf{size}($A$); s=\tbf{length}($\mbf{blocks}$); $high=\mbf{blocks}(1)$; \\
$[Q,R]$=\tbf{local\_qr}($A(:,1:{high})$);\\
\tbf{for} \= $k=2:s$ \\
\>$low= high+1$; $high=high+\mbf{blocks}(k)$;\\
\>$[Q_{new},R_{new},S_{new}]$=\tbf{block\_CGS2\_step}($Q,A(:,low :  high)$);\\
\>$R = \left(\begin{array}{cc} R & S_{new} \\ 0 & R_{new} \end{array}\right)$; $Q = \left(\begin{array}{cc} Q & Q_{new} \end{array}\right)$; \\
\tbf{end}; \\
\tbf{end}; {\em block\_CGS2} \\
\end{tabbing}

\section{Error analysis of the algorithm BCGS2}
 \label{sec:CGS2erranalysis}
  Our error analysis of Function 2.5  is the result of the error analysis of one step of Function 2.4
followed by an induction argument. The details of our proof use standard error analysis assumptions and techniques.

\subsection{Error Bounds for Function 2.4}

To establish our error bound for Function 2.4, we establish two bounds.

The first, on
\begin{equation}
\normtwo{B-U S_B -Q_B R_B},
\label{eq:resid}
\end{equation}
has no preconditions. The second, on
\begin{equation}
\normtwo{U^T Q_B}
\label{eq:orth}
\end{equation}
requires one of two assumptions.
 The first assumption is
\begin{equation}
0< \macheps f_{sing}(m,t,p) \normtwo{B} \normtwo{R_B^{-1}} \leq \gamma
\label{eq:assumption1}
\end{equation}
where
\begin{eqnarray}
\gamma &=& L_F(m,t,p)/f_1(m,t,p) <1, \label{eq:gammadef} \\
f_{sing}(m,t,p) &=& f_1(m,t,p)+L_F(m,t,p)+\gamma L_5(p). \label{eq:fsing}
\end{eqnarray}
The second is that
\begin{equation}
\normtwo{R_2^{-1}} \leq \ysqrt{1+\gamma^2}
\label{eq:assumption2}
\end{equation}
for $\gamma$ in (\ref{eq:gammadef}).

The first theorem covers (\ref{eq:resid}).

\begin{theorem}
\label{thm:resid}
Assume the hypothesis and notation of Lemma \ref{lem:BCGS_step}. Then the computed $Q_B$,$R_B$ and
$S_B$ from Function 2.4
satifies
\begin{equation}
\normtwo{B-US_B -Q_R R_B} \leq \macheps f_{resid}(m,t,p)\normtwo{B} \hot
\label{eq:residbnd}
\end{equation}
where
\begin{equation}
f_{resid}(m,t,p)= 2 L_1(m,p)+2 L_3(t,p)+L_4(p)+L_5(p).
\label{eq:fresid}
\end{equation}
\end{theorem}

{\bf Proof.} Using Lemma \ref{lem:BCGS_step_backward} on the second block CGS step in Function 2.4
yields
\begin{eqnarray*}
Q_B R_2 &=& (I_m-\XXT{U}) Q_1 +F_2\\
&=& Q_1 -US_2 -U(\delta S_2) + F_2 \\
&=& Q_1 -US_2+\Delta Y_2-\delta Y_2.
\end{eqnarray*}
Multiplying by $R_2$ yields
\[
Q_B R_2 R_1 = (I_m-\XXT{U}) Q_1 R_1 + F_2 R_1
\]
thus
\[
Q_B R_B = (I_m-\XXT{U}) Q_1 R_1 + F_2 R_1 - Q_B (\delta R_B).
\]
If we use the fact that
$F_2 =\Delta Y_2 -\delta Y_2 +U(\delta S_2)$ from Lemma \ref{lem:BCGS_step_backward}, then we have
\[
Q_B R_B = Q_1 R_1 -US_2 R_1  + (\Delta Y_2 -\delta Y_2) R_1- Q_B (\delta R_B).
\]
Expanding $Q_1 R_1$ using Lemma \ref{lem:BCGS_step_backward} yields
\begin{eqnarray*}
Q_B R_B &=& (I_m-\XXT{U})B -U S_2 R_1 +F_1  + (\Delta Y_2 -\delta Y_2) R_1- Q_B (\delta R_B)\\
&=& B-US_1 - US_2 R_1 -U(\delta S_1) +F_1 + (\Delta Y_2 -\delta Y_2) R_1- Q_B (\delta R_B).
\end{eqnarray*}
Using the definition of $F_1$ in Lemma \ref{lem:BCGS_step_backward} and the backward error for $S_B$ in (\ref{eq:deltaSB}) we
have
\[
Q_B R_B = B-US_B - U(\delta S_B) +\Delta Y_1 -\delta Y_1   + (\Delta Y_2 -\delta Y_2) R_1- Q_B (\delta R_B).
\]
Using norm bounds yields
\begin{eqnarray*}
\normtwo{B-US_B -Q_B R_B} &\leq& \normtwo{\delta S_B} + \normtwo{\Delta Y_1}+\normtwo{\delta Y_1} + +\normtwo{\delta Y_2}\normtwo{R_1} +
\normtwo{\delta Y_2} \normtwo{R_1}+\normtwo{\delta R_B}\\
&\leq& \macheps( [L_4(p)+L_1(m,p)+L_3(t,p)+L_5(p)]\normtwo{B} +[L_1(m,p) +L_3(t,p)]\normtwo{R_1}) \hot \\
&\leq& \macheps(2 L_1(m,p)+2 L_3(t,p)+L_4(p)+L_5(p))\normtwo{B} \hot\\
&=& \macheps f_{resid}(m,t,p)\normtwo{B} \hot
\end{eqnarray*}
establishing (\ref{eq:residbnd})--(\ref{eq:fresid}).

\medskip

A crucial norm relationship is given by the following lemma.

\begin{lemma}
\label{lem:R2UQ1}
Let $R_2$ and $Q_1$ be result of implementing Function 2.4
in floating point arithmetic with machine unit $\macheps$.
Then, if $R_2$ is nonsingular, for $\gamma <1$
\begin{equation}
\normtwo{U^T Q_1 R_2^{-1}} \leq \gamma +\mathcal{O}(\macheps)
\label{eq:UQ1R2rel1}
\end{equation}
if and only if
\begin{equation}
\normtwo{R_2^{-1}} \leq (1+\gamma^2)^{1/2} +\mathcal{O}(\macheps).
\label{eq:UQ1R2rel2}
\end{equation}
\end{lemma}

{\bf Proof.}
We start with interpreting Lemma \ref{lem:BCGS_step} for the
second CGS step in in Function 2.4
which leads to
\[
Q_B R_2 = (I_m-\XXT{U})Q_1  + F_2.
\]
Taking the normal equations matrices of both sides yields
\begin{equation}
R_2^T \normal{Q_B} R_2 = Q_1^T (I_m-\XXT{U})^2 Q_1 + F_2^T(I_m-\XXT{U})Q_1+Q_1^T (I_m-\XXT{U}) F_2 + \normal{F_2}.
\label{eq:R2normal}
\end{equation}
An expansion of $Q_1^T (I_m-\XXT{U})^2 Q_1$ produces
\begin{eqnarray}
Q_1^T (I_m-\XXT{U})^2 Q_1 &=& \normal{Q_1}-Q_1^T UU^T Q_1-Q_1^T U(I_t-\normal{U})U^T Q_1 \nonumber \\
&=& I-Q_1^T UU^T Q_1\\
&+&\normal{Q_1}-I-Q_1^T U(I_t-\normal{U})U^T Q_1 \label{eq:Q1Uexpand}
\end{eqnarray}
so that the combination of (\ref{eq:R2normal}) and (\ref{eq:Q1Uexpand}) is
\begin{equation}
\normal{R_2} = I-Q_1^T UU^T Q_1 + E
\label{eq:R2Q1E}
\end{equation}
where
\begin{eqnarray*}
E&=& E_1+E_2+E_3, \\
E_1 &=&  F_2^T(I_m-\XXT{U})Q_1+Q_1^T (I_m-\XXT{U}) F_2 + \normal{F_2},\\
E_2 &=& \normal{Q_1}-I_p-Q_1^T U(I_t-\normal{U})U^T Q_1, \\
E_3 &=& R_2^T (\normal{Q_B}-I_p) R_2.
\end{eqnarray*}
Since
\begin{eqnarray*}
\normtwo{E_1} &\leq& 2\normtwo{ F_2} + \normtwo{F_2}^2 \leq 2\macheps L_F(m,t,p) \hot, \\
\normtwo{E_2} &\leq & \normtwo{I_p-\normal{Q_1}}+\normtwo{I_t-\normal{U}} \hot \leq \macheps [L_1(m,p)+f_1(m,t,p)]\hot,\\
\normtwo{E_3} &=& \normtwo{I_p-\normal{Q_B}} \normtwo{R_2}^2 \leq \macheps L_1(m,p) \hot,
\end{eqnarray*}
we have
\begin{eqnarray*}
\normtwo{E} &\leq& \normtwo{E_1}+\normtwo{E_2} + \normtwo{E_3} \\
&\leq& \macheps [f_1(m,t,p)+2L_F(m,t,p)+2L_1(m,p)]\hot.
\end{eqnarray*}

Now to show the equivalence between (\ref{eq:UQ1R2rel1}) and (\ref{eq:UQ1R2rel2}).
Since we assume that $R_2$ is nonsingular, we can rewrite (\ref{eq:R2Q1E}) as
\[
I_p = R_2^{-T} R_2^{-1} -  R_2^{-T}Q_1 U U^T Q_1 R_2^{-T}+ R_2^{-T}E R_2^{-1}
\]
so that
\begin{equation}
 R_2^{-T} R_2^{-1}=I_p+  R_2^{-T}Q_1 U U^T Q_1 R_2^{-T}- R_2^{-T}E R_2^{-1}.
\label{eq:R2Q1R2}
\end{equation}
If $\lambda_1(\cdot)$ is the leading eigenvalue of the contents, then
(\ref{eq:R2Q1R2}) is given by
\begin{equation}
\lambda_1( R_2^{-T} R_2^{-1})=1+  \lambda_1(R_2^{-T}Q_1 U U^T Q_1 R_2^{-1})+  \xi \normtwo{R_2^{-1}}^2
\label{eq:eigenbnd}
\end{equation}
where
\begin{equation}
|\xi| \leq \normtwo{E} .
\label{eq:xibnd}
\end{equation}
Using the relationship, $\lambda_1(\normal{C})=\normtwo{C}^2$ on (\ref{eq:eigenbnd}) yields
\begin{equation}
\normtwo{R_2^{-1}}^2(1-\xi) =1+ \normtwo{U^T Q_1 R_2^{-1}}^2 .
\label{eq:normtwobnd}
\end{equation}
Thus assuming (\ref{eq:UQ1R2rel1}) yields
\[
\normtwo{R_2^{-1}}^2= [1+(\gamma+\mathcal{O}(\macheps))^2]/(1-\xi)= 1+\gamma^2 +\mathcal{O}(\macheps)
\]
establishing (\ref{eq:UQ1R2rel2}). Likewise, a similar algebraic manipulation of  (\ref{eq:normtwobnd}) shows that (\ref{eq:UQ1R2rel2}) implies (\ref{eq:UQ1R2rel1}).

\medskip

Before showing the effect of assumption  (\ref{eq:assumption1}), we need a small
technical lemma.

\begin{lemma}
\label{lem:R1R2}
Assume  \eqref{eq:assumption1}  and assume the hypothesis and notation of Lemma 2.1.
 Then Function 2.4
 produces nonsingular $R_1$ and $R_2$.
\end{lemma}

{\bf Proof.}
Since the smallest singular value of $R_B$ satisfies
\[
\sigma_p(R_B) = \normtwo{R_B^{-1}}^{-1}
\]
assumption (\ref{eq:assumption1}) may be written
\[
\gamma \sigma_P(R_B) \geq \macheps f_{sing}(m,t,p)\normtwo{B} > 0.
\]
>From (\ref{eq:singvalbnd}),
\[
\sigma_p(R_B) -  L_5(p) \macheps \leq \sigma_p(R_2 R_1)
\]
thus
\begin{equation}
0 <\macheps [f_1(m,t,p) +L_F(m,t,p)]\normtwo{B} \leq \gamma  \sigma_p(R_2 R_1).
\label{eq:sR2R1}
\end{equation}
Thus $R_2 R_1$ is nonsingular. Since $R_1$ and $R_2$ are square, $R_1$ and $R_2$
are each nonsingular.

\medskip

We now show the effect of the assumption (\ref{eq:assumption1}).

\begin{lemma}
\label{lem:UQ1R2}
Assume  \eqref{eq:assumption1}  and assume the hypothesis and notation of Lemma \ref{lem:BCGS_step}.
 Then Function 2.4
 produces $Q_1$ and $R_2$ such that \eqref{eq:UQ1R2rel1}
and \eqref{eq:UQ1R2rel2} hold.
\end{lemma}

{\bf Proof.} From Lemma \ref{lem:R2UQ1}, (\ref{eq:UQ1R2rel1}) and (\ref{eq:UQ1R2rel2}) are equivalent  so we need only prove (\ref{eq:UQ1R2rel1}).

  Interpreting Lemma \ref{lem:BCGS_step} and using the result of Lemma \ref{lem:R1R2} for the first CGS step in Function 2.4
yields
\begin{eqnarray*}
U^T Q_1 R_2^{-1} &=& U^T (I_m-\XXT{U}) B R_2^{-1}R_1^{-1} + U^T F_1 R_1^{-1} R_2^{-1}\\
&=& [(I_t-\normal{U})U^T B +U^T F_1](R_2 R_1)^{-1}.
\end{eqnarray*}
Thus
\begin{eqnarray}
\normtwo{U^T Q_1 R_2^{-1}} &\leq& [\normtwo{I_t-\normal{U}}\normtwo{U^T B} + \normtwo{U^T F_1}]\normtwo{(R_2 R_1)^{-1}} \nonumber \\
&\leq& \macheps [f_1(m,t,p) +L_F(m,t,p)]\normtwo{B}\normtwo{(R_2 R_1)^{-1}} \hot. \label{eq:normBR2R1}
\end{eqnarray}

Combining (\ref{eq:sR2R1}) and (\ref{eq:normBR2R1}) yields
\[
\normtwo{U^T Q_1 R_2^{-1}} \leq \gamma  \sigma_p(R_2 R_1) \normtwo{(R_2 R_1)^{-1}} \hot =\gamma\hot
\]
which satisfies (\ref{eq:UQ1R2rel1}).

\medskip

Now we have a conditional bound on $\normtwo{U^T Q_B}$ from Function \ref{fun:BCGS2_step}.

\begin{theorem}
\label{thm:UQB}
Assume the hypothesis and notation of Lemma \ref{lem:BCGS_step}. Assume also that
$U$ satisfies \eqref{eq:Uorth} and either \eqref{eq:assumption1} or \eqref{eq:assumption2}
holds. Then Function 2.4
produces $Q_B$ that satisfies
\begin{equation}
\normtwo{U^T Q_B} \leq \macheps [1+ \sqrt{2}] L_F(m,t,p) \hot.
\label{eq:UTBbnd}
\end{equation}
\end{theorem}

{\bf Proof.}
Applying  Lemma \ref{lem:BCGS_step_backward} to the second \tbf{block\_CGS\_step} of
Function 2.4
yields
\[
Q_B R_2 = (I_m-\XXT{U}) Q_1 + F_2, \quad \normtwo{F_2} \leq \macheps L_F(m,t,p) \hot.
\]
Thus
\[
U^T Q_B =(I_t-\normal{U})U^T Q_1 R_2^{-1} + U^T F_2 R_2^{-1}.
\]
Norm bounds lead to
\[
\normtwo{U^T Q_B} \leq \normtwo{I_t-\normal{U}}\normtwo{U^T Q_1 R_2^{-1}} + \normtwo{ F_2}\normtwo{ R_2^{-1}}.
\]
From Lemmas and \ref{lem:R2UQ1} and  \ref{lem:UQ1R2},
either assumption (\ref{eq:assumption1}) or (\ref{eq:assumption2}) yields (\ref{eq:UQ1R2rel1})--(\ref{eq:UQ1R2rel2}), thus
\begin{equation}
\normtwo{U^T Q_B} \leq \gamma\normtwo{I_t-\normal{U}}  + (1+\gamma^2)^{1/2} \normtwo{ F_2} \hot.
\label{eq:UTQBbnd1}
\end{equation}
which from (\ref{eq:gammadef}) becomes
\begin{eqnarray}
\normtwo{U^T Q_B} &\leq& \macheps [\gamma f_1(m,t,p)  + (1+\gamma^2)^{1/2} L_F(m,t,p)] \hot \nonumber \\
&\leq& \macheps [(1+(1+\gamma^2)^{1/2}) L_F(m,t,p)] \hot \label{eq:UTQBbnd2}
\end{eqnarray}
Using $\gamma <1$ produces (\ref{eq:UTBbnd}).

\subsection{Error Bounds for Function 2.5}
\label{sec:blockCGS2}

Obtaining the bounds (\ref{eq:IQTQbnd}) and (\ref{eq:AQRbnd}) are simply the result of
induction arguments on Theorems \ref{thm:resid} and \ref{thm:UQB}.

In the arguments of this section, we assume that all of the blocks $A_1,\ldots, A_s$ have the same dimension, i.e., $p_1 = \cdots = p_s=p$. To have blocks
of differing size, we could just assume that $p=\max_{1 \leq i \leq s} p_i$ and
make some other minor adjustments to the proofs in this section.

 We begin with (\ref{eq:AQRbnd}) and let $t_k = kp$.

\begin{theorem}
\label{thm:AQRbnd}
Assume the hypothesis and notation of Lemma \ref{lem:BCGS_step}. Let $\hat{Q}_k$ in \eqref{eq:hatQk} and
$R_k$ in  \eqref{eq:Rk} be the result of $k$ steps of Function 2.5.
Then for $\numlist{k}{1}{s}$,  $A_k$ in \eqref{eq:hatQk} satisfies
\begin{equation}
\normtwo{A_k -\hat{Q}_k R_k} \leq \macheps f_2(m,t_{k-1},p) \normtwo{A_k} \hot
\label{eq:Akresid}
\end{equation}
where
\begin{equation}
f_2(m,t_{k-1},p)= k^{1/2} f_{resid}(m,t_{k-1},p).
\label{eq:f2spec}
\end{equation}
Thus \eqref{eq:AQRbnd} follows from \eqref{eq:Akresid}--\eqref{eq:f2spec} by taking $k=s$.
\end{theorem}

{\bf Proof.}
For $k=1$, $\hat{A}_1 = A_1$ $\hat{Q}_1 = Q_1$ and $R_1=R_{11}$, by our assumption (\ref{eq:BQRbnd})
\begin{eqnarray*}
\normtwo{\hat{A}_1 - \hat{Q}_1 R_1 } &\leq& \macheps L_1(m,p)\normtwo{A_1} \hot\\
&\leq& \macheps  f_{resid} (m,0,p) \normtwo{A} \hot \\
&\leq& \macheps f_{2}(m,0,p) \normtwo{A} \hot.
\end{eqnarray*}

For the induction, step assume that the theorem holds up to step $k$ and prove
for step $k+1$. For $k <s$, we have that
\[
\hat{A}_{k+1} - \hat{Q}_{k+1} (k+1)p = \left(\begin{array}{cc} \hat{A}_k & A_{k+1} \end{array}\right) - \left(\begin{array}{cc} \hat{Q}_k & Q_{k+1} \end{array}\right) \left(\begin{array}{cc} R_k & S_{k+1} \\ 0 & R_{k+1,k+1} \end{array}\right)
\]
where $S_{k+1} \in \mathbb{R}^{t_k \times p}$ and $R_{k+1,k+1} \in \mathbb{R}^{p \times p}$.
Thus,
\begin{equation}
\normtwo{\hat{A}_{k+1} - \hat{Q}_{k+1} R_{k+1}}^2 \leq \normtwo{\hat{A}_{k} - \hat{Q}_{k} R_{k}}^2
+ \normtwo{A_{k+1} - \hat{Q}_k S_{k+1} - Q_{k+1} R_{k+1,k+1} }^2.
\label{eq:residrecc}
\end{equation}
The first term is bounded by the induction hypothesis, the second results from applying Theorem \ref{thm:resid} to $A_{k+1}$, $\hat{Q}_k$, $S_{k+1}$, $Q_{k+1}$ and $R_{k+1,k+1}$ which gives us
\begin{eqnarray*}
\normtwo{\hat{A}_{k} - \hat{Q}_{k} R_{k}} &\leq&\macheps k^{1/2} f_{resid}(m,t_{k-1},p) \normtwo{\hat{A}_k }\hot \label{eq:term1bnd} \\
\normtwo{A_{k+1} - \hat{Q}_k S_{k+1} - Q_{k+1} R_{k+1,k+1} } &\leq&\macheps f_{resid}(m,t_k,p) \normtwo{A_k}\hot
\label{eq:term2bnd}
\end{eqnarray*}

Combining (\ref{eq:residrecc}), (\ref{eq:term1bnd}), and (\ref{eq:term2bnd}) yields
\begin{equation}
\normtwo{\hat{A}_{k+1} - \hat{Q}_{k+1} R_{k+1}}^2 \leq \macheps^2 [k f_{resid}^2(m,t_{k-1},p)\normtwo{\hat{A}_k}^2 +f_{resid}^2(m,t_k,p)\normtwo{A_{k+1}}^2]+\mathcal{O}(\macheps^3) .
\label{eq:bothtermsbnd}
\end{equation}

Since $f_{resid}(\cdot)$ is monotone nondecreasing is all of its arguments and $\normtwo{\hat{A}_k}$,
$\normtwo{A_{k+1}} \leq \normtwo{\hat{A}_{k+1}}$,
 (\ref{eq:bothtermsbnd}) becomes
\[
\normtwo{\hat{A}_{k+1} - \hat{Q}_{k+1} R_{k+1}}^2 \leq (k+1)\macheps^2 f_{resid}^2(m,t_k,p) \normtwo{\hat{A}_{k+1}}^2+\mathcal{O}(\macheps^3) .
\]
Taking square roots establishes the induction step of the argument.

To prove the orthogonality bound (\ref{eq:IQTQbnd}), we need to make define
\begin{eqnarray}
\gamma_k &\stackrel{def}{=}& \frac{L_F(m,t_{k-1},p)}{f_1(m,t_{k-1},p)}, \quad \numlist{k}{2}{n}
\label{eq:gammakdef1} \\
&=& \ysqrtn{\alpha^2(k-1)+1}, \quad\alpha =\sqrt{7+4\sqrt{2}} \approx 3.56.
\label{eq:gammakdef2}
\end{eqnarray}
Let $R_2^{(k)}$ be the upper triangular matrix produce in the second call to Function 2.3
in the $kth$ step of Function 2.5. Our generalizations of assumptions (\ref{eq:assumption1}) and (\ref{eq:assumption2}) to Function 2.5 are
\begin{eqnarray}
f_{sing}(m,t_{k-1},p) \normtwo{A_k} \normtwo{R_{kk}^{-1}} \leq \gamma_k
\label{eq:assumption1gen} \\
\normtwo{[R_2^{(k)}]^{-1}} \leq \ysqrt{1+\gamma_k^2}, \numlist{k}{2}{n}
\label{eq:assumption2gen}
\end{eqnarray}
where $f_{sing}(m,t,p)$ is defined by (\ref{eq:fsing}). Using these two assumptions,
we have our final theorem.

\begin{theorem}
\label{thm:IQTQbnd}
Assume the hypothesis and notation of Theorem \ref{thm:AQRbnd} and that
either assumption \eqref{eq:assumption1gen} or \eqref{eq:assumption2gen} holds.
Then for $\numlist{k}{1}{n}$,
\begin{equation}
\normtwo{I_{t_k} -\normal{\hat{Q}_k}} \leq \macheps f_1(m,t_{k-1},p) \hot
\label{eq:IQTQbndk}
\end{equation}
where $f_1(\cdot)$,$\gamma_k$ and $\alpha$ satisfy \eqref{eq:gammakdef1}--\eqref{eq:gammakdef2}.
Interpreting \eqref{eq:IQTQbndk} for $k=s$ yields \eqref{eq:IQTQbnd}.
\end{theorem}

{\bf Proof.} This is a proof by induction on Theorem \ref{thm:UQB}.
For $k=1$, we note that $\hat{Q}_1 = Q_1$ which just results from
$[Q_1,R_{11}]=\mbf{local\_qr}(A_1)$. Thus
\begin{eqnarray*}
\normtwo{I_{p}-\normal{\hat{Q}_1}} & =& \normtwo{I_{p}-\normal{Q_1}} \\
&\leq& \macheps L_1(m,p) \hot \\
&\leq& \macheps f_1(m,0,p) \hot
\end{eqnarray*}

For the induction step, assume  $t_k \leq n$  that (\ref{eq:IQTQbndk})
holds for $k$. Then
\[
I_{(k+1)p}-\normal{\hat{Q}_{k+1}} = \left(\begin{array}{cc} I_{t_k} - \normal{\hat{Q}_k} & \hat{Q}_k^T Q_{k+1} \\
Q_{k+1}^T \hat{Q}_k & I_{p} - \normal{Q_{k+1}} \end{array}\right).
\]
so that
\begin{eqnarray*}
\normtwo{I_{(k+1)p}-\normal{\hat{Q}_{k+1}}}&=& \normtwo{\left(\begin{array}{cc} I_{t_k} - \normal{\hat{Q}_k} & \hat{Q}_k^T Q_{k+1} \\
Q_{k+1}^T \hat{Q}_k & I_{p} - \normal{Q_{k+1}} \end{array}\right)} \\
&\leq& \normtwo{\left(\begin{array}{cc} \normtwo{I_{t_k} - \normal{\hat{Q}_k}} & \normtwo{\hat{Q}_k^T Q_{k+1}} \\
\normtwo{Q_{k+1}^T \hat{Q}_k} & \normtwo{I_{p} - \normal{Q_{k+1}}} \end{array}\right)}
\end{eqnarray*}
Invoking the induction hypothesis, applying Theorem \ref{thm:UQB} to $\hat{Q}_k$ and $Q_{k+1}$,
and using the assumption (\ref{eq:IQBbnd}) yields
\begin{eqnarray}
\normtwo{I_{(k+1)p}-\normal{\hat{Q}_{k+1}}}&\leq&  \macheps \normtwo{ \left(\begin{array}{cc} f_1(m,t_k,p) & (1+\sqrt{2})L_F(m,t_k,p)\\(1+\sqrt{2})L_F(m,t_k,p) & L_1(m,p) \end{array}\right)}   \hot
\nonumber \\
&\leq& \macheps \normF{ \left(\begin{array}{cc} f_1(m,t_{k-1},p) & (1+\sqrt{2})L_F(m,t_k,p) \\ (1+\sqrt{2})L_F(m,t_k,p) & L_1(m,p)\end{array}\right) }\hot \nonumber \\
&\leq& \macheps c_{orth}(m,t_k,p) \hot \label{eq:f1exp}
\end{eqnarray}
where using the implicit definition of $f_1(\cdot)$  the definitions of $\gamma_k$ and $\alpha$  in (\ref{eq:gammakdef1})--(\ref{eq:gammakdef2}), and of $L_F(\cdot)$ in (\ref{eq:Fbnd2}), we have
\[
c_{orth}(m,t_k,p) \stackrel{def}{=}\ysqrt{ f_1^2(m,t_{k-1},p)+ 2((1+\sqrt{2})^2L_F^2(m,t_k,p)
+L_1^2(m,p)}
\]
We can bound $c_{orth}(\cdot)$ by
\begin{eqnarray*}
c_{orth}(m,t_k,p) &\leq& \ysqrt{ f_1^2(m,t_{k-1},p)+2((1+\sqrt{2})^2L_F^2(m,t_k,p)
+L_F^2(m,t_k,p)}\\
&\leq& \ysqrt{ f_1^2(m,t_{k-1},p)+(7+4\sqrt{2})L_F^2(m,t_k,p)}\\
&=& \ysqrt{ f_1^2(m,t_{k-1},p)+ \alpha^2 L_F^2(m,t_k,p)}\\
&=& \ysqrt{\gamma_k^{-2} L_F^2(m,t_{k-1},p)+\alpha^2 L_F^2(m,t_k,p)}.
\end{eqnarray*}
Since $L_F(\cdot)$ is nondecreasing in all of its arguments
\begin{eqnarray}
c_{orth}(m,t_k,p) &\leq& \ysqrt{\gamma_k^{-2} + \alpha^2}L_F(m,t_k,p)
\nonumber \\
&=& \gamma_{k+1}^{-1} L_F(m,t_k,p) =f_1(m,t_k,p)
\label{eq:inductionf1}
\end{eqnarray}
Combining (\ref{eq:f1exp})--(\ref{eq:inductionf1}) yields the induction step for
(\ref{eq:IQTQbndk}).

\subsection{Interpreting the Bounds for Function 2.5}
\label{sec:CGS2}

 Function 2.2  is just Function 2.5 with
$p=1$ and with \tbf{local\_qr} producing  $R_B = \left( r_b \right)$ and
$Q_B = \left( \mbf{q}_b \right)$ from $B=\left( \mbf{b} \right)$ from the normalization
\[
r_b =\| \mbf{b} \|_2; \quad \mbf{q}_b = \mbf{b}/r_b.
\]
When $p=1$, $t_k =k$ and we will interpret it as such.

In floating point arithmetic, the computed values satisfy
\begin{equation}
r_b =\| \mbf{b} \|_2 (1+\delta), \quad |\delta | \leq (m/2+1) \macheps \hot
\label{eq:rberr}
\end{equation}
and
\begin{equation}
\mbf{q}_b = (I_m+ D) \mbf{b}/r_b, \quad D=diag(d_i), \quad \|D\|_2 \leq \macheps.
\label{eq:qberr}
\end{equation}

From (\ref{eq:rberr})--(\ref{eq:qberr}), it follows that
\[
\|\mbf{b} - \mbf{q}_b r_b\|_2 \leq \macheps \|\mbf{b}\|_2
\]
and
\[
|1-{\mbf{q}_b}^T {\mbf{q}_b} | \leq (m+4)\macheps \hot
\]
which are  (\ref{eq:IQBbnd})--(\ref{eq:BQRbnd}) with $L_1(m,1) =\max\{ m+4,1\} =m+4$.

The other operations of a CGS step are
\[
\bar{\mbf{s}}+ \delta \bar{\mbf{s}} = U^T \mbf{b}
\]
where
\[
\normtwo{\delta \bar{\mbf{s}}}\leq\macheps mk^{1/2} \normtwo{\mbf{b}} \hot.
\]
Thus, $L_2(m,k,1)= mk^{1/2}$.

 We also have
\[
\bar{\mbf{y}}+ \delta \bar{\mbf{y}} = \mbf{b} -U\bar{\mbf{s}}
\]
where
\[
\normtwo{\delta \bar{\mbf{y}}} \leq \macheps (1+k^{3/2})\normtwo{\mbf{b}} \hot.
\]
Thus $L_3(k,1)=1+k^{3/2}$.
Thus one step is
\[
\bar{\mbf{q}} \bar{r} = (I_m-\XXT{U}) \mbf{b} + \mbf{f}
\]
where
\[
\normtwo{\mbf{f}} \leq \macheps L_F(m,k,1) \normtwo{\mbf{b}} \hot
\]
and
\[
L_F(m,k,1)= L_1(m,1)+L_2(m,k,1)+L_3(k,1)= m+k^{1/2}(m+k)+4.
\]
The two operations
\begin{eqnarray*}
\mbf{s}_B+\delta \mbf{s}_B &=& \mbf{s}_1 + \mbf{s}_2 r_1 \\
r_B + \delta r_B &=& r_2 r_1
\end{eqnarray*}
satisfy
\begin{eqnarray*}
\normtwo{\delta \mbf{s}_B} &\leq& 2 \macheps \normtwo{\mbf{b}} \hot \\
|\delta r_B| &\leq& \macheps \normtwo{\mbf{b}}\hot
\end{eqnarray*}
so that $L_4(1)=2$ and $L_5(1)=1$. Thus
\[
f_{resid}(m,k,1)= 2(m+4)+2(1+k^{3/2})+3=2m+2k^{3/2}+13.
\]

Taking $t=s=n$, we have that
 Function 2.2 obtains a Q--R factorization
satisfying
\begin{eqnarray*}
\normtwo{A-QR} &\leq& \macheps  n^{1/2} f_{resid}(m,n,1) \normtwo{A} \hot \\
&\leq& \macheps(2mn^{1/2} +2n^2 +13n^{1/2})\normtwo{A} \hot.
\end{eqnarray*}

The condition for near orthogonality of $Q$ has a nice interpretation.
Assumption (\ref{eq:assumption1gen}) may be written
\begin{equation}
\macheps f_{sing}(m,k,1)\normtwo{\mbf{a}_k} \leq \gamma_k |r_{kk}|
\label{eq:CGS2assumption1}
\end{equation}
where again,
\[
 \gamma_k = L_F(m,k,1)/f_1(m,k,1) = \ysqrtn{\alpha^2 (k-1)+1}
\]
 and
\begin{eqnarray*}
f_{sing}(m,k,1)&=& f_1(m,k,1)+L_F(m,k,1)+L_5(1) \\
                &=& [\ysqrt{\alpha^2(k-1)+1}+1]L_F(m,k,1) +2  \\
                &=& [\ysqrt{\alpha^2(k-1)+1}+1] ( m+k^{1/2}(m+k)+4)+2.
\end{eqnarray*}
Assumption (\ref{eq:assumption2gen}) is
\begin{equation}
|r_{kk}^{(2)}| \geq 1/\ysqrt{1+\gamma_k^2}.
\label{eq:CGS2assumption2}
\end{equation}
where $r_{kk}^{(2)}$ is the diagonal element in the second Gram--Schmidt step at step
$k$ of Function 2.2.

Either assumption leads to  the bound
\[
\normtwo{I-\normal{Q}} \leq \macheps f_1(m,n,1) \hot
\]
where
\[
f_1(m,n,1) =\gamma_n^{-1} L_F(m,n,1).
\]

Notice that (\ref{eq:CGS2assumption1}) is merely an assumption that each of  the diagonals of $R$
is sufficiently bounded away from $\macheps \normtwo{\mbf{a}_k}$. There is no assumption on the condition number
of $R$ (or $A$) and (\ref{eq:CGS2assumption1}) much weaker  than the assumption given by Giraud et al. \cite{Gir05}
for Function 2.2. The second assumption, (\ref{eq:CGS2assumption2}), is very
similar to an assumption discussed by Abdelmalek \cite{Abd71}.

\section{Conclusions}
\label{sec:conclusions}

Function 2.5 is a new block classical Gram--Schmidt Q--R factorization with
reorthogonalization. We have shown that as long as the diagonal blocks on $R$ do not become too
ill--conditioned, the factorization produces a near orthogonal $Q$ according to the criterion
(\ref{eq:IQTQbnd}) and a small residual according to the criterion (\ref{eq:AQRbnd}).

Moreover, if we consider the block size $1$, we have improved a bound of Giraud et al. \cite{Gir05}
for Function 2.2 showing that a near left orthogonal $Q$ is produced if diagonals of
$R$ are bounded sufficiently away from zero.


\end{document}